\documentclass[11pt,bezier]{article}
\usepackage{amsmath,amssymb,amsfonts,euscript}
\usepackage[usenames]{color}
\usepackage{graphicx}

\newpage

\newcommand{\be}{\begin{equation}}
\newcommand{\ee}{\end{equation}}
\newtheorem{theorem}{Theorem}[section]
\newtheorem{proof}{Proof}[section]

\newtheorem{remark}{Remark}[section]

\title{\bf\Large   Markov Switching  asymmetric  GARCH Model: Stability and Forecasting}
\vspace{-1cm}
\author
{N. Alemohammad\thanks{\scriptsize Department of Mathematics and Computer Science, Shahed University, Tehran , Iran.}, S. Rezakhah\thanks{\scriptsize   Faculty of Mathematics and Computer Science, Amirkabir University of Technology,
Tehran, Iran.} \thanks{\scriptsize Corresponding author, Email: rezakhah@aut.ac.ir}, S.H. Alizadeh\thanks{ \scriptsize  Department of Computer Engineering and IT, Qazvin
Branch, Islamic Azad University, Qazvin, Iran.
}}

\begin{document}
\maketitle
\begin{abstract}

\end{abstract}
\begin{abstract}
	A new Markov switching asymmetric GARCH  model is proposed  where each state follows a logistic smooth transition  structure between effects of  positive and negative shocks.  This consideration provides better forecasts in many financial time series. The asymptotic finiteness of the second moment is investigated. The parameters of the model are estimated by applying MCMC methods  through Gibbs  and griddy Gibbs sampling.  Applying the log return of some part of 
  S$\&$P 500 indices, we show the    competing performance of in sample fit and out of sample forecast
 volatility and value at risk 
of the proposed model.  The Diebol-Mariano test shows that  the presented   model  outperforms all competing models in forecast volatility.

	 \quad\\
	 
	 
\end{abstract}
 \section{Introduction}
 Volatility modeling in financial time series has been widely studied over past few decades. The ARCH  and
GARCH models , introduced by Engle \cite{engle a} and Bollerslev \cite{bollerslev}, are surely the most popular classes of volatility models. Hamilton and Susmel \cite{hamilton} introduced the Markov-Switching GARCH (MS-GARCH) by merging GARCH model with a hidden Markov chain, where each state allows a different GARCH behavior.  Such structure  improves  forecasting of volatilities. Gray \cite{gray}, Klaassen \cite{klaassen} and Haas, et al. \cite{haas a} proposed some different variants of MS-GARCH models.    For further studies on MS-GARCH models, see  Abramson and Cohen \cite{abramson}, Ardia \cite{ardia b}, Alemohammad et. al \cite{alemohammad}  and Bauwens et al. \cite{bauwens c}.  \\

\par
{ One restriction  of the GARCH model is its  symmetry to the sign of past shocks. This is improved by letting the conditional variance to be a function of size and sign  of the preceding observation.
 Such studies  started by  Black \cite{black}, who investigated  asymmetric effects of  positive and negative shocks on volatilities.
  This consideration is important in financial markets as there exists higher volatility in response to bad news (negative shocks). \cite{gonzalez-rivera}.
  Study of the  asymmetric GARCH  started by  Engle \cite{engle b} and continued as the Exponential GARCH (EGARCH) model by Nelson \cite{nelson},  GJR-GARCH model  by Glosten, et al.\cite{glosten} and Threshold GARCH (TGARCH) model by Zakoian \cite{zakoian}.
  The other   asymmetric structures are Smooth transition models introduced by   Gonzalez-Rivera \cite{gonzalez-rivera},  {Ardia \cite{ardia b},} Medeiros and Veiga \cite{medeiros},  and Haas et al. \cite{haas b}.}

   \par
\textcolor{blue}{In this paper, we study some   Markov switching  GARCH model where  the  volatility in  each regime is coupled with the smooth transition between the effects of negative and positive schocks. 
 More precisely,   the presented model  
	  considers  different smooth transition structure by states  where   each state   describes some  time dependent convex combination  between  asymmetric  effects of positive and   negative and  shocks. 
The new model obviates the absence of asymmetric property in the Markov switching GARCH model and switching 
	between  different levels of volatility in the Smooth transition GARCH model  (ST-GARCH),   presented by   Lubrano \cite{lubrano}. }

\textcolor{blue}{ Ardia \cite{ardia b} considered some  MS-GARCH model  where   the  asymmetric  effects of  volatilities   are considered   by applying 
  indicator functions of some  predefined  non-positive thresholds in states. This cause a sudden shift of volatility structure at corresponding threshold in each state.       
    Alemohammad et. al \cite{alemohammad} considered a Markov  {switching} GARCH model  where the smooth transition between structures for high and low  volatilities   are  in effect of   size  of the preceding  return.
 In this paper we study the case where the volatility structure   in  each state  follows some  smooth transition between the effects 
 of  positive and negative shocks  based on the  
 preceding log return. So it is expected to provide much better fitting, especially  
when smooth transitions between such effects are evident. }
	As such model employs  all past observations, we reduce the volume of calculations  by proposing a dynamic programming algorithm.  We also derive  sufficient condition for stability of the model by applying 
  the method of Abramson and Cohen \cite{abramson} and Medeiros \cite{medeiros}.
	\\ \textcolor{blue}{ Markov chain Monte Carlo (MCMC) methods are widely applied for parameter estimation of the regime switching GARCH models, see, e.g, \cite{bauwens c} and \cite{ardia b}. As nonlinearity is considered, the likelihood function becomes tricky to maximize since it is hard to differentiate, see  \cite{lubrano}. In addition the existence of latent variables 
also makes  the Bayesian method to be required.   The  privileges of  the MCMC  methods is  that  avoid the common problem of local maxima encountered in the Maximum likelihood (ML) estimation, see \cite{ardia a} (section 7.7). 
   The parameters of our model are estimated by applying MCMC methods through  Gibbs and Griddy-Gibbs sampling. 
We present  a simulation example to show the competitive performance of our model in compare to GARCH,  MS-GARCH and  ST-GARCH models.  Using  $S\&P$500 indices from  $3/01/2005$ to $3/11/2014$   we show that our model has much better fitting by providing less forecasting error,  better performance base on Diebold-Mariano test and  value at risk of out-of-sample forecasting  of one day ahead volatility, 
in compare with  GARCH, MS-GARCH,  EGARCH,GJR-GARCH and ST-GARCH models. 
 We also show that  our model outperforms the competing models for in-sample fit by using the Deviance information criterion.}

 \par The Markov switching smooth transition GARCH model is presented in section 2. Section 3 is devoted to the statistical properties of the model. Estimation of the parameters of the model are studied in section 4. \textcolor{blue}{ Simulation studies are followed to show  competing performance  of presented model  in section 5}. Section 6 is dedicated to the analysis of the efficiency of the proposed model by applying the model to the S$\&$P 500  indices  for  $3/01/2005$ to $3/11/2014$ . Section 7 concludes.

\section{  Markov switching asymmetric  GARCH model}

We consider the Markov   switching  smooth transition GARCH model, in summary  MS-STGARCH as 
 \begin{equation}\label{1}
 y_t=\varepsilon_{t}\sqrt{H_{Z_t,t}},\hspace{2cm}
\end{equation}
where $\{\varepsilon_{t}\}$ are iid  standard normal variables, $\{Z_t\}$ is an irreducible and aperiodic  Markov chain on finite state space $E=\{1,2,\cdots,K\}$
 with transition probability matrix
$\; P=||p_{ij}||_{K\times K},\;$
where  $\, p_{ij}=p(Z_t=j|Z_{t-1}=i),\; i,j \in \{1,\cdots,K\}$,
and stationary probability measure
$\,
\pi=(\pi_1,\cdots,\pi_K)^{\prime}.$   Also given that $Z_t=j$, $H_{j,t}$ (the conditional variance of regime j)  is defined as
\begin{equation}\label{2}
H_{j,t}=a_{0j}+y^2_{t-1}(d_{j,t})+b_jH_{j,t-1},
\end{equation}

where
\begin{equation}\label{3}
d_{j,t}=a_{1j}(1-w_{j,t-1})+a_{2j}w_{j,t-1}
\end{equation}
and 
 the weights ($w_{j,t}$) are  logistic function of the past observation as
\begin{equation}\label{4}
w_{j,t-1}=\frac{1}{1+\exp(-\gamma_jy_{t-1})}\ \ \ \  \gamma_j>0,\ j=1,\cdots,K,\ \ \ \
\end{equation}

 \noindent  which  are monotonically increasing with respect to previous observation and are bounded , $0<w_{j,t-1}<1$. The parameter $\gamma_j>0$ is called the slope parameter. The weight function   $w_{j,t-1}$ goes to one when $y_{t-1}\rightarrow +\infty$ and so $d_{j,t}$ tends to $a_{2j}$. Also it  goes to zero when $y_{t-1}\rightarrow -\infty$ and so $d_{j,t}$ tends to $a_{1j}$.  Therefore  the effect of negative shocks are mainly described   by $a_{1j}$  and of positive shocks by  $a_{2j},j=1, \cdots, K$.
  As often   negative shocks have greater effect on volatilities  than positive ones,  one could assume that $a_{1j}>a_{2j}$ in each regime. To impose the idea in model building, we recommend to consider a higher prior for  $a_{1j}$ in each state and any estimation procedure.
  Logistic weight functions in  states have the potential to describe different speed for smooth transitions which are in effect of $\gamma_j,j=1,\cdots,K$ and also different effect limits as $a_{1j}$ and $a_{2j}$. This enables one to provide a flexible model for describing such different transitions.   Plots of such  logistic  weight functions  for
 the returns of  S$\&$P 500 indices, which are studied later in this paper, are presented in Figure 1.
 Indeed in each regime the coefficient of $y^2_{t-1}$ is time dependent that causes the volatility structure being under the influence size and sign of the observations and it makes distinct from GARCH model.

    As $\gamma_j\rightarrow \infty$, the logistic weight function considers a step function for positive and negative shocks. When $\gamma_j$ approaches to zero, $w_{j,t-1}$ goes  to $1/2$ and the  MS-STGARCH model tends to the Markov switching GARCH model(MS-GARCH). In the case of single regime,  our model is the  smooth transition GARCH (STGARCH) model that is introduced by  Lubrano \cite{lubrano}.\\
 \par It is assumed that $\{\varepsilon_{t}\}$ and  $\{Z_t\}$ are independent.  Sufficient conditions to guarantee strictly positive conditional variance (\ref{2})  are that
 $a_{0j}$ to be positive and $a_{1j},a_{2j},b_{j}$
being nonnegative.

\label{fig1}\input{epsf}
\epsfxsize=3in \epsfysize=2in
 \begin{figure}
\vspace{0in}

\centerline{\hspace{0in}\epsfxsize=6.4in \epsfysize=2.8in  \epsffile{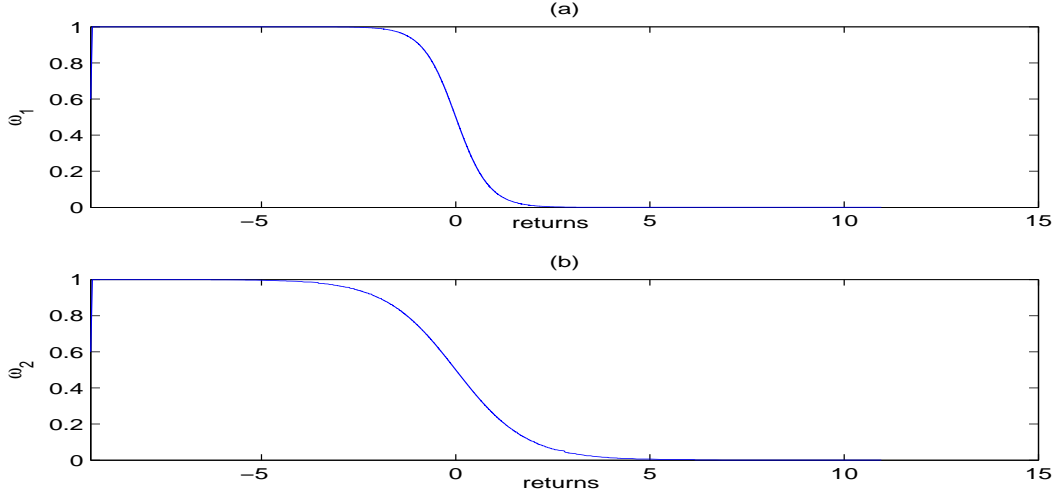}}
\vspace{0in}

\caption{\scriptsize Logistic function ($w_{j,t}$) for $S\&P$ returns  (a): the weight of first regime,    (b): the weight of second regime}
\end{figure}

\section{Statistical Properties of the model}
In this section,  the statistical properties of the  MS-STGARCH model are investigated and   the conditional density and  variance  of the process is obtained. As the evaluation of the asymptotic behavior of the  second moment  in our model isn't so easy to follow,  we apply the method of Abramson and Cohen \cite{abramson} and Medeiros \cite{medeiros} to obtain an appropriate upper bound for the asymptotic value of the unconditional variance to show its stability.
\subsection{Conditional density and variance}
Let  $\mathcal{I}_{t}$  be the information  up to time t. Following the method of Alemohammad et al.  \cite{alemohammad}, the conditional density function of $y_t$ given past information can be written as\\
\begin{equation}\label{7}
	f(y_t|\mathcal{I}_{t-1})= \sum_{j=1}^{K}{\alpha_j^{(t)}\phi(\frac{y_t}{\sqrt{H_{j,t}}})}
\end{equation}
where    $\phi(.)$ is the probability density function of the standard normal distribution and  $\alpha_j^{(t)}=p(Z_t=j|\mathcal{I}_{t-1})$, that is obtained in the following remark.
{{{\remark{{\label{9}}
				The value of $\alpha_j^{(t)}$  is obtained recursively by
				\begin{equation}
				\alpha_j^{(t)}=\frac{\sum_{m=1}^{K}{f(y_{t-1}|Z_{t-1}=m,\mathcal{I}_{t-2})}\alpha_m^{(t-1)}p_{m,j}}{\sum_{m=1}^{K}{f(y_{t-1}|Z_{t-1}=m,\mathcal{I}_{t-2})}\alpha_m^{(t-1)}},
				\end{equation}
where $p_{mj}=p(Z_t=j|Z_{t-1}=m), m,j=1,\cdots,K$ are the transition probabilities and
	\begin{equation*}
f(y_{t-1}|Z_{t-1}=m,\mathcal{I}_{t-2})=\phi(\frac{y_{t-1}}{\sqrt{H_{m,t-1}}}).
\end{equation*}		}}
			{{\proof  See Appendix A. }}
\\
\par The conditional variance of the MS-STGARCH model is given by
\begin{align}\label{8}
Var(Y_t|\mathcal{I}_{t-1}) &=\sum_{j=1}^{K}{\alpha_j^{(t)}H_{j,t}}
\end{align}
as $H_{j,t}$ is the conditional variance of j-th state. This relation shows that the conditional variance of this model is affected by changes in regime and conditional variance of each state. \\As using all past observations for forecasting could increase the complexity of the model, we reduce the volume of calculations  by proposing a dynamic programming algorithm.
At each time $t$, $\alpha_i^{(t)}$ (in equation (\ref{7}), (\ref{8})) can be obtained from a dynamic programming method based on the forward recursion algorithm, proposed in remark (\ref{9}).
\subsection{Stability}
To show  the asymptotic wide sense  stationarity of the MS-GARCH model, Abramson and Cohen \cite{abramson}  evaluated unconditional second moment  by conditioning on past observations and  providing a linear equation between present and  past volatilities. Then they showed that the necessary and sufficient condition for the asymptotic wide sense stationarity  of the model is that 
the spectral radius of some related matrix to be less than one.  For the MS-STGARCH, this method fails as the logistic weights cause that  the evaluation of  unconditional variance of observations to be so complicated.  Lubrano \cite{lubrano}     obtained a recurrsive relation for the volatility of STGARCH model that by which and by imposing some condition 
they showed  the stationarity and persistence of the volatility of  STGARCH Model.  Surely while we have different structure for the volatilities of states in MS-STGARCH  their method is not applicable. 
 So we study the asymptotic boundedness of the second order moment, which shows the stability of the model. 
In this subsection, we investigate the stability of second moment of the  MS-STGARCH model. So it would be enough to find   an upper bound of  the second moment of the process, see \cite{alemohammad}. Let M be a positive constant and 
 \begin{equation}\label{16}
 {\bf{\Omega}}=[a_{01}+|a_{21}-a_{11}|M^2,\cdots,\\a_{0K}+|a_{2K}-a_{1K}|M^2)]^{\prime},
 \end{equation}
  be a vector with K component,  ${{C}}$ denotes  a $K^2$-by-$K^2$ block matrix as \\
  \begin{equation}\label{17}
  {{C}}=\left(
             \begin{array}{cccc}
              {{C}_{11}}& {{C}_{21}}& \cdots & {{C}_{K1}} \\
               {{C}_{12}} & {{C}_{22}} & \cdots & {{C}_{K2}} \\
               \vdots &  &  & \vdots\\
               {{C}_{1K}} & {{C}_{2K}} & \cdots & {{C}_{KK}} \\
             \end{array}
           \right)
  \end{equation}
  where
 \begin{equation}\label{18}
 {{C}_{jk}}=p(Z_{t-1}=j|Z_t=k)( {{u}}{{e}}^{\prime}_j+{{v}}),\ \ \ \ \ \ \quad\ j,k=1,\cdots,K,
 \end{equation}
 $ {{u}}=[a_{11}+(\delta+\frac{1}{2})|a_{21}-a_{11}|,\cdots,a_{1K}+(\delta+\frac{1}{2})|a_{2K}-a_{1K}|]^{\prime}$,   ${ {e}}_j$ is a K-by-1 vector that the jth element of it is one and other components are zero, and $v=\parallel v_{ij}\parallel_{i,j=1}^K$ is a diagonal  matrix that $v_{jj}=b_j$ for $j=1,\cdots,K$.\\
Let $\Pi=[\pi_1{ {e}}^{\prime}_1,\cdots,\pi_K{ {e}}^{\prime}_K]$ }and   $\rho(A)$ denotes  that  spectral radius of  matrix A. Now we present the following theorem regarding the stability condition of the MS-STGARCH model. 

{\theorem Let $\{Y_t\}_{t=0}^{\infty}$ follows the MS-STGARCH  model, defined by (\ref{1})-(\ref{4}),
 the process is asymptotically stable in  variance  and $\lim_{t\rightarrow\infty}E(Y^2_t)\leq{ {\Pi^\prime(I-{ {C}})^{-1}\dot{{ {\Omega}}}}}$, if $\rho({ {C}})
<1.$}\\
{{\proof  See Appendix B. }}

\section{Estimation}
For the estimation of parameters,  we apply the Bayesian MCMC method, that is extensively used in literature (\cite{bauwens b}, \cite{bauwens c} and \cite{lubrano}).
\par Let $Y_t=(y_1,\cdots,y_t)$ and $Z_t=(z_1,\cdots,z_t)$ be the samples of observations and hidden variables respectively.  We consider two states for the model  with parameters  $\theta=(\theta_1,\theta_2)$, where $\theta_k=(a_{0k},a_{1k},a_{2k},b_{k},\gamma_k)$ for $k=1,2$ and transition probabilities
 $\eta=(\eta_{11},\eta_{12}, \eta_{21}, \eta_{22})$ where $\eta_{ij}=p(z_{t+1}=j|z_t=i)$. The posterior density can be represented as

\begin{equation}\label{21}
p(\theta,\eta,Z|Y)\propto p(\theta,\eta)p(Z|\theta,\eta)f(Y|\theta,\eta,Z),
\end{equation}
 where $Y=(y_1,\cdots,y_T)$, $Z=(z_1,\cdots,z_T)$, T is the total number of samples and $p(\theta,\eta)$ is the prior density. By assuming that the value of $z_1$ is known, conditional probability mass function of $Z$ given the $(\theta,\eta)$ is independent of $\theta$, so
\begin{align}
p(Z|\theta,\eta)=& p(Z|\eta_{11},\eta_{22})
\nonumber\\
&=\prod_{t=1}^{T}{p(z_{t+1}|z_t,\eta_{11},\eta_{22})}
\nonumber\\
&=\eta_{11}^{n_{11}}(1-\eta_{11})^{n_{12}}\eta_{22}^{n_{22}}(1-\eta_{22})^{n_{21}},
\end{align}
where $n_{ij}=\#\{z_t=j|z_{t-1}=i\}$ (the number of transitions from regime i to regime j).
 The conditional density function of $Y$ given the realization of $Z$ and the parameters is factorized in the following way:
\begin{equation}
f(Y|\eta,\theta,Z)=\prod_{t=1}^{T}{f(y_t|\theta,z_t=k,Y_{t-1})},\ \ \ k=1,2,
\end{equation}
where the one step ahead predictive densities are:
\begin{equation}\label{22}
f(y_t|\theta,z_t=k,Y_{t-1})=\frac{1}{\sqrt{ 2\pi H_{k,t}}}\exp(-\frac{y^2_t}{2H_{k,t}}).
\end{equation}
\par Since the straight sampling from the posterior density (\ref{21}) is not possible, we apply the Gibbs sampling algorithm for three blocks: $\theta$, $\eta$ and $Z$. \\
In implementing Gibbs algorithm, we consider  the superscript $(r)$ on a parameter to denote its value at the r-th iteration of the algorithm. At any iteration of the algorithm, three steps are considered:\\
(i) Draw the random sample of the state variable $Z^{(r)}$  given  $,\eta^{(r-1)},\ \theta^{(r-1)}$.\\
(ii) Draw the random sample of the transition probabilities $\eta^{(r)}$ given $Z^{(r)}$.\\
(iii) Draw the random sample of the $\theta^{(r)}$ given $Z^{(r)}$ and $\eta^{(r)}$.\\

These steps are repeated until the convergency is obtained. In what follows the sampling of each block are explained. \\
\subsection{Sampling $z_t$}
This step is devoted to the sampling of the conditional probability $p(z_t|\eta,\theta, Y_t)$ which considered  by Chib\cite{chib}, see also \cite{kaufman}.  Suppose $p(z_1|\eta,\theta, Y_0,)$ be the stationary distribution of the chain, then\\
\begin{equation}
p(z_t|\eta,\theta, Y_t)\propto f(y_t|\theta,z_t=k,Y_{t-1})p(z_t|\eta,\theta, Y_{t-1}),
\end{equation}
where the predictive density $ f(y_t|\theta,z_t=k,Y_{t-1})$ is calculated by (\ref{22}) and by the law of total probability, $p(z_t|\eta,\theta, Y_{t-1})$ is given by\\
\begin{equation}
p(z_t|\eta,\theta, Y_{t-1})=\sum_{z_{t-1}=1}^{K}{p(z_{t-1}|\eta,\theta, Y_{t-1})\eta_{z_{t-1}z_t}},
\end{equation}
where K is the number of states. Given the  probabilities ($p(z_t|\eta,\theta, Y_t)$), we run a backward algorithm, starting from $t=T$,  $z_T$ is derived from $p(z_T|\eta,\theta,Y)$. For $t=T-1,\cdots,0$ the corresponding samples are derived from $p(z_t|z_{t+1},\cdots,z_T,\theta,\eta,Y)$,which satisfies 
\begin{equation*}
p(z_t|z_{t+1},\cdots,z_T,\theta,\eta,Y)\propto p(z_t|\eta,\theta,Y_t)\eta_{z_t,z_{t+1}}.
\end{equation*}
Derive  $z_t$ from $p(z_t|.)=p_{z_t}$  by the following procedure:\\
	first evaluate $q_{j}=p(Z_t=j|Z_t\geq j,.)$ by $$p(Z_t=j|Z_t\geq j,.)=\frac{p_j}{\sum_{l=j}^{K}{p_l}},$$
	then  generate a number u from the standard uniform distribution (U(0,1)). If $u\leq q_j$ then put $z_t=j$ otherwise increasing $j$ to $j+1$ and generate another u from U(0,1) and repeat this step by comparing this with $q_{j+1}$. 
\subsection{Sampling $\eta$}
This stage is devoted to sample $\eta=(\eta_{11},\eta_{22})$ from the posterior probability $p(\eta|\theta, Y_t, Z_t)$ that is independent of $Y_t, \theta$. We consider independent beta prior density for each of   $\eta_{11}$ and $\eta_{22}$.  So,
$$p(\eta_{11}|Z_t)\propto p(\eta_{11})p(Z_t|\eta_{11})=\eta_{11}^{c_{11}+n_{11}-1}(1-\eta_{11})^{c_{12}+n_{12}-1},$$
where $c_{11}$ and $c_{12}$ are the parameters of beta prior, $n_{ij}$ is the number of transition from $z_{t-1}=i$ to $z_t=j$.  In the same way  the sample of $\eta_{22}$ is obtained.
\subsection{Sampling $\theta$}
The posterior density of $\theta$ given the prior $p(\theta)$ is given by:
\begin{equation}\label{23}
p(\theta|Y, Z)\propto p(\theta) \prod_{t=1}^{T}{f(y_t|\theta,z_t=k,Y_{t-1})}= p(\theta)\prod_{t=1}^{T}{\frac{1}{\sqrt{ 2\pi H_{k,t}}}\exp(-\frac{y^2_t}{2H_{k,t}})},
\end{equation}
which is independent of $\eta$. To sample from the $p(\theta|Y, Z)$ we use the Griddy Gibbs algorithm that introduced by Ritter and Tanner \cite{ritter}. This method has had wide application in literature, see \cite{bauwens a} , \cite{bauwens b} and \cite{bauwens c}.
 Given samples at iteration $r$ the Griddy Gibbs at iteration $r+1$ proceeds as follows:\\
\\
1. Select a grid of points, such as $a_{0i}^1,a_{0i}^2,\cdots,a_{0i}^G$. Use (\ref{23}) to evaluate the kernel of conditional posterior density function 
of $a_{0i}$ given all the values of Z, Y and $\theta$ except $a_{0i}$ $k(a_{0i|Z_t,Y_t,\theta_{-a_{0i}}})$ over the grid points to obtain the vector $G_k=(k_1,\cdots,k_G)$.\\
2. By a deterministic integration rule using the G points, compute
$G_{\Phi}=(0,\Phi_2,\cdots,\Phi_G)$
with
\begin{equation}
\Phi_j=\int_{a_{0i}^1}^{a_{0i}^j}{k(a_{01}|\theta_{-a_{0i}}^{(r)},Z^{(r)}_t,Y_t)da_{0i}},\ \ \ i=2,\cdots,G.
\end{equation}
3. Simulate $u\sim U(0,\Phi_G)$ and invert $\Phi(a_{0i}|\theta_{-a_{0i}}^{(r)},Z^{(r)}_t,Y_t)$ by numerical interpolation to obtain a sample $a_{0i}^{(r+1)}$ from $p(a_{0i}|\theta_{-a_{0i}}^{(r)},Z^{(r)}_t,Y_t)$.\\
4. Repeat steps 1-3 for other parameters.
\\
\par Prior densities of  elements of $\theta$ can be considered  as independent uniform densities over finite intervals.
\label{fig2}\input{epsf}
\epsfxsize=3in \epsfysize=2in
\begin{figure}
	\centerline{\epsfxsize=5.7in \epsfysize=1.8in \epsffile{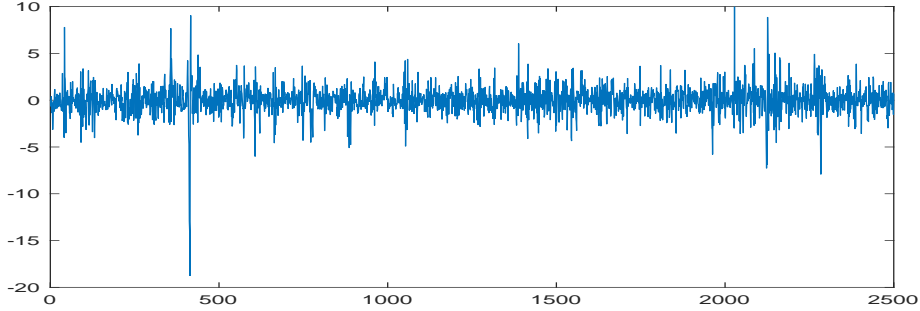}}
	\vspace{-0.25in}
	\caption{\scriptsize  Simulated data.}
\end{figure}
\textcolor{blue}{\section{Simulation results}
We have simulated 2500 sample from the proposed model   (\ref{1})-(\ref{4}) for two states, $j=1,2$. Figure 2 shows  the  plot of the simulated time  series and Table 1  reports some descriptive  statistics of the simulated data. \\
Using the Bayesian inference, we estimate the parameters of the MS-STGARCH  by applying  the first  2000 samples.  The prior density of each parameter is assumed to be uniform  over a finite interval except for transition probabilities  $\eta_{11}$ and $\eta_{22}$ which  are drawn from some  beta distribution. Table 2 demonstrates the true  values of the parameters and also posterior means and standard deviations of the corresponding estimators over 10000 iterations, which 5000 of them are discarded as burn-in samples.}  \textcolor{blue}{ The results of this table shows that the mean square errors (MSE) of the estimated parameters are adequately small.  Using simulated data, we compare the in-sample fit and out-of-sample forecasting performance of  the presented  model with GARCH, STGARCH and MS-GARCH models in subsetions 5.1 and 5.2 respectively. 
 }


\begin{table}
	\vspace{.05in}
	\caption{{\small Descriptive statistics  for the simulated data.}}
	
	\hspace{-2cm}{\hspace{.8in}\small
		\begin{tabular}{c c c c c c c}
			\hline
			& Mean & Std. dev. & Skewness & Maximum & Minimum &  Kurtosis  \\
			
			\hline
		Simulated data&-0.03 &1.470 &-1.621 &9.057  &-18.754  &24.623 \\
			
			\hline
		\end{tabular}}
	\end{table}

	\begin{table}
		\vspace{-.05in}
		\caption{{\scriptsize Results of the Bayesian Estimation of the  simulated MS-STGARCH model. }}
		\hspace{4cm}{\scriptsize
			\begin{tabular}{c c c c c}
				\hline
				
				& True values &  Mean & Std. dev. & MSE \\
				\hline
				$a_{01}$ &0.300 &0.314  & 0.031 & 0.0001 \\
				$a_{11}$ &0.200 &0.224  & 0.022 & 0.0005 \\
				$a_{21}$ &0.050  &0.12  & 0.014& 0.0002\\
				$b_{1}$ &0.500  &0.510  & 0.049 & 0.002  \\
			$\gamma_1$ &1.500 &1.517 & 0.179 & 0.032\\
				$a_{02}$ & 1.900 &1.750  &0.09 & 0.008 \\
				$a_{12}$ &0.700  &0.619  & 0.060 & 0.004\\
				$a_{22}$ &0.100 &0.094  & 0.016 & 0.0003\\
				$b_{2}$ & 0.250 &0.217 &0.021 & 0.0004  \\
							  
				$\gamma_2$ &0.500 & 0.619 & 0.013 & 0.0002\\
				$\eta_{11}$ & 0.970& 0.982 & 0.093 & 0.009 \\
				$\eta_{22}$ &0.850  &0.853  &0.097 & 0.009 \\
				\hline
				
			\end{tabular}}
		\end{table}
\subsection{In sample performance analysis}
\begin{table}[!hbp]
	
	\caption{{\small Deviance information criterion (DIC) for simulated data}}
	\hspace{3.5cm}\scriptsize
	\begin{tabular}{cc}
		\hline
		{\small Model} &{\scriptsize{ \small S$\&$P 500 returns }}\\
		\hline
		{\small GARCH} & 5896\\
		
		{\small ST-GARCH} & 6021.4\\
		{\small MS-GARCH} &5887.8 \\
		{\small MS-STGARCH} & 5858.6*\\
		\hline
	\end{tabular}
\end{table}
\textcolor{blue}{In order to compare the goodness of fit of the proposed model with GARCH, STGARCH and MS-GARCH, we apply the deviance information criterion (DIC) introduced by Spiegelhalter et al \cite{speigelhalter}. DIC is a Bayesian version  of the reputable Akaike information criterion (AIC) that  designed specifically for Bayesian estimation
	that involves MCMC simulations, see \cite{gelman} and \cite{stata} . The smallest DIC determines the best model. Berg applied DIC for the familiy of stochastic volatility (SV) models \cite{berg} and Ardia for the family of asymmetric GARCH models  \cite{ardia b}.  } \textcolor{blue}{ 
	\\ In the Markov-switching models, the likelihood is calculated by the following formula:
	\begin{equation*}
	f(Y|\Theta)=\prod_{t=1}^{T}{f(y_t|\mathcal{I}_{t-1},\Theta)}
	\end{equation*}
	in which $\Theta$ is the vector of all parameters in model and $f(y_t|\mathcal{I}_{t-1})$ is obtained from (\ref{7}). The deviance information criterion  is computed as:
	\begin{equation}
	DIC=2 \log(f(Y|\hat{\Theta}))-4E_{Y|\Theta}[\log(f(Y|\Theta))],
	\end{equation}
	that $\hat{\Theta}$ is the posterior means of the vector $\Theta$. The results concerning DIC of the simulated data are reported in Table 3. It is apparent that the DIC of MS-STGARCH is the smallest value in the table. Thus our considered model has the best fit to the simulated data set among competing models.}
\subsection{Out-of-sample forecasting performance analysis}
For appraising the performance of MS-STGARCH in forecasting, we survey the one-day-ahead value at risk (VaR) forecasts for the last 500 data of the  simulated data. The one-day-ahead value at risk level $\alpha\in(0,1)$, VaR$(\alpha)$ is obtained by calculating the $(1-\alpha)$th
percentile of the one-day-ahead predictive distribution (4.14).
To test the VaR at level $\alpha$, we evaluate the sequence $\{V_t(\alpha)\}$ by
\begin{displaymath}
V_t(\alpha)=\left\{ \begin{array}{ll}
I\{y_{t+1}<VaR (\alpha)\} & \textrm{if  $\alpha>0.5$}\\
I\{y_{t+1}> VaR (\alpha)\} & \textrm{if $\alpha\leq 0.5$}.
\end{array} \right.
\end{displaymath}
The out-of-sample VaR  at level $\alpha$ has good performance if the sequence  $\{V_t(\alpha)\}$ are  independent and  obey the following  distribution


\begin{displaymath}
V_t(\alpha)\sim\left\{ \begin{array}{ll}
Bernoulli(1-\alpha) & \textrm{if  $\alpha>0.5$}\\
Bernoulli(\alpha) & \textrm{if $\alpha\leq 0.5$},
\end{array} \right.
\end{displaymath}

The three likelihood ratio statistics for  unconditional coverage ($LR_{uc}$), independence ($LR_{ind}$) and conditional coverage ($LR_{cc}$)  are as follows \cite{christofferssen}:\\
1.  LR statistic for the test of unconditional coverage,
\begin{equation*}
LR_{uc}=-2\ln[\frac{\phi^{n_1}(1-\phi)^{n_0}}{\hat{\pi}^{n_1}(1-\hat{\pi})^{n_0}}]\sim\chi^2_{(1)},
\end{equation*}
where $\phi$ is the parameter of related Bernoulli distribution, which could be $1-\alpha$ or $\alpha$, $n_1$ is the number of 1's and $n_0$ is the number of 0's in the $V_t(\alpha)$ series and $\hat{\pi}=\frac{n_1}{n_1+n_0}$.\\
2. LR statistic for the test of independence,
\begin{equation*}
LR_{ind}=-2\ln[\frac{\hat{\pi}_{*}^{n_{00}+n_{10}}(1-\hat{\pi}_{*})^{n_{11}+n_{01}}}{\hat{\pi}_{1}^{n_{00}}(1-\hat{\pi}_{1})^{n_{01}}\hat{\pi}_{2}^{n_{11}}(1-\hat{\pi}_{2})^{n_{10}}}]\sim\chi^2_{(1)},
\end{equation*}
where $n_{ij}$ is the number of transition from i to j ($i,j=0,1$) in the $V_t(\alpha)$ series, $\hat{\pi}_{1}=\frac{n_{00}}{n_{00}+n_{01}}$, $\hat{\pi}_{2}=\frac{n_{11}}{n_{10}+n_{11}}$ and $\hat{\pi}_{*}=\frac{n_{00}+n_{10}}{n_{00}+n_{01}+n_{10}+n_{11}}$.\\
\\
3.LR statistic for the test of conditional coverage,

$$LR_{cc}=LR_{ind}+LR_{uc},$$
$LR_{cc}$ has $\chi^2$ distribution with two degrees of freedom.  When the value of $LR_{cc}$ is less than the critical value of $\chi^2$ distribution one infer  that the  conditional coverage is correct and there exist good VaR forecasts. \\
 The results of the tests for simulation example  are reported in Table 4. The second and third  columns demonstrate  the theoretical expected violations and the number of empirical violations respectively. 
 \textcolor{blue}{According to the results of Table 4, at the $ 5\%$  significance levels, as  $\chi^2_{1,0.95}=3.841$, the  $LR_{uc}$ test is
 	rejected four times  for GARCH and  STGARCH, two times for MS-GARCH   and one time for MS-STGARCH models.
 	For some risk levels the test of independence (IND test) is not applicable since no consecutive violations have been occurred. In such  cases  $n_{00}=0$ and so the  $LR_{ind}$ statistic  becomes infinity.  The $LR_{ind}$ statistic at $ 5\%$ significance level is bigger than critical  value for one case of GARCH and MS-STGARCH and also two cases of STGARCH.   The conditional coverage (CC) test is higher than  critical value $\chi^2_{2,0.95}=5.991$ with two degrees of freedom three times for the GARCH, four times for STGARCH, two times for the MS-GARCH and one time for the MS-STGARCH.} 
 \begin{table}
 	\vspace{-.05in}
 	\caption{{\scriptsize VaR results  of  simulated data.}}
 	\hspace{-.5cm}{\hspace{.8in}\scriptsize
 		\begin{tabular}{c c c c c c c}
 			\hline
 			Model   & $\alpha$ & $E(V_t(\alpha))$ & N & UC & IND &  CC  \\

 			\hline
 			&0.99 &5 &13& 8.973  & 4.243   & 13.216 \\
 			&0.95 &25 & 42 & 10.195  & 1.756 & 11.95\\
 			&0.9 &50 & 62 & 2.997 & 0.03 & 3.027\\
 			GARCH   &0.1& 50 & 50 & 0 & 0.2664&  0.2664\\
 			&0.05& 25& 37 & 5.317 & 0.256& 5.573\\
 			&0.01& 5& 12& 7.111& 1.15 & 8.263 \\
 			\hline
 			&0.99&5& 12& 7.111& 4.854 & 11.964 \\
 			&0.95&25& 38& 6.181& 1.527 & 7.706 \\
 			&0.9&50& 63& 3.499 & 0.005& 3.504\\
 			STGARCH&0.1&50& 53& 0.197& 0.030& 0.227\\
 			&0.05&25& 39& 7.102& 0 & 7.103 \\
 			&0.01&5&12& 7.111& 4.854 & 11.964 \\
 			\hline
 			&0.99&5& 13& 8.970 & 0.914 & 9.887 \\
 			&0.95&25&  34& 3.080& 2.81& 5.892\\
 			&0.9&50& 56& 0.773& 0.558& 1.330\\
 			MS-GARCH&0.1&50& 50& 0 & 0.266&  0.266\\
 			&0.05&25& 38 & 6.181 & 0.357 & 6.538\\
 			&0.01&5&10&  3.914& 1.75 & 5.665 \\
 			\hline
 			&0.99 &5 & 8 & 1.538  & NA  & NA \\
 			
 			&0.95 &25 & 36 & 4.510& 3.987 & 8.498\\
 			&0.9 &50 &62& 2.996 & 0.015 & 3.011\\
 			MS-STGARCH   &0.1&50& 54& 0.347& 0.274& 0.622\\
 			&0.05&25& 34& 3.081& 0.215 & 3.295\\
 			&0.01&5& 9& 2.613& 2.126 & 4.739 \\
 			\hline
 			
 		\end{tabular}}
 	\end{table}

\begin{table}
\vspace{.05in}
\caption{{\small Descriptive statistics  for the  S$\&$P 500  index daily log returns.}}

\hspace{-1cm}{\hspace{.8in}\small
\begin{tabular}{c c c c c c c}
  \hline
   & Mean & Std. dev. & Skewness & Maximum & Minimum &  Kurtosis  \\

  \hline
 S$\&$P 500&0.023 &1.287 &-0.337 &10.957  &-9.469  &14.049 \\

      \hline
  \end{tabular}}
\end{table}
\section{Empirical data set}
  By applying   daily log returns of the S$\&$P 500 for the period of  03/01/2005 to 03/11/2014 (2500 observations), we compare the performance of our model with the GARCH, MS-GARCH and ST-GARCH ones. From the 2500 observations of S$\&$P 500, the first 2000 observations are employed to estimate the parameters and the remaining 500
samples are used for forecasting analysis.  Figure 3  plots the daily log returns in percentages\footnote{log return in percentage is
	defined as $r_t=100*\log(\frac{P_t}{P_{t-1}})$, where $P_t$ is the index level at time t.}  of the S$\&$P500  indices. 
 
   In Table 5,   the descriptive statistics of  the log returns in percentages
 	  are presented.
This table shows that the means are close to zero and there are some slightly negative skewness  and  excess kurtosis for the data set.
\label{fig2}\input{epsf}
\epsfxsize=3in \epsfysize=2in
 \begin{figure}
\centerline{\epsfxsize=5.7in \epsfysize=1.8in \epsffile{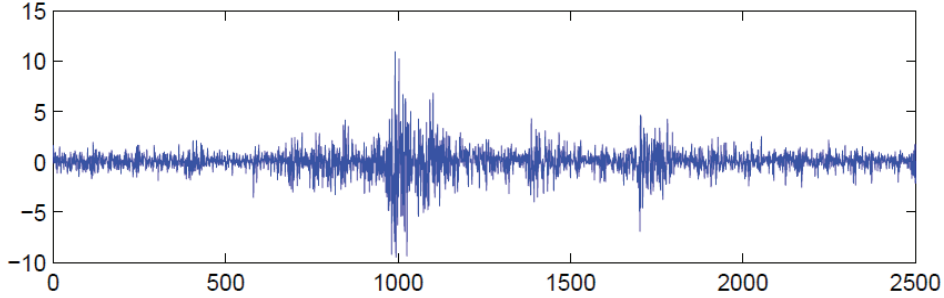}}
\vspace{-0.25in}
\caption{\scriptsize   Percentage daily log returns of  S$\&$P 500 data.}
\end{figure}
    \begin{table}
\vspace{.2in}
\caption{{\small Posterior means and standard deviations (S$\&$P 500  daily log returns). }}

\hspace{-0.7cm}{\hspace{.4in} \small
\begin{tabular}{lllllllllllllllll}
\hline
  \multicolumn{2}{c}{}$\hspace{0.5cm}$MS-STGARCH $\hspace{-1cm}$  & \multicolumn{2}{c}{}$\hspace{-.8cm}$MS-GARCH $\hspace{-1.2cm}$  & \multicolumn{2}{c}{}$\hspace{1cm}$ST-GARCH $\hspace{-1cm}$  & \multicolumn{2}{c}{}$\hspace{1cm}$GARCH\\
\cline{2-2}\cline{4 -4}\cline{6-6}\cline{8-8}
 &Mean  &Std.dev.   &Mean &Std.dev & Mean & Std.dev  &Mean &Std.dev\\
\hline
 $a_{01}$& 0.194 &0.001    &0.233 &0.004  & 0.336&0.011 & 0.269&0.005 &\\
 $a_{11}$& .276 &0.008     &0.278  & 0.012   &0.421 &0.019 & 0.120& 0.007\\
 $a_{21}$&0.085  & 0.006   &0  &0   &0.121 &0.016  &0&0\\
 $b_{1}$& 0.289 &0.003    &0.320 &0.009   &0.364 &0.014 & 0.439&0.004\\
  $\gamma_1$&2.345  &0.132   &0  &0  &2.206 &0.218 & 0&0 \\
 $a_{02}$  &0.717  &0.087   &0.779 & 0.012  &- &- &- &-\\
 $a_{12}$&0.677 &0.007   &0.430 & 0.011   &- &- &  -&-\\
 $a_{22}$&0.365  & 0.013    &0 &0   &- &-  & -&- \\
 $b_{2}$&0.264  & 0.015    &0.207 &0.007   &- &-  & -&-\\
 $\gamma_2$&1.097  & 0.017   &0 &0  &- &-  & -&- \\
 $\eta_{11}$& 0.986 & 0.004  &0.994 & 0.003   &- &-  & -&-\\
 $\eta_{22}$& 0.985 &0.005   & 0.991 & 0.004   &- &-  & -&-\\
 \hline
\end{tabular}}
\end{table}
\begin{table}
		\vspace{-.05in}
		\caption{{\scriptsize Parameter estimation of the EGARCH and GJR-GARCH models. }}
		\hspace{4cm}{\scriptsize
			\begin{tabular}{c c c }
				\hline
				
			Coefficients	& EGARCH & GJR-GARCH   \\
				\hline
				Constant &0.007 &0.019   \\
				arch effect  &0.127 &0   \\
				garch effect &0.975  &0.903  \\
				leverage effect &-0.147  &0.162   \\
				\hline

			\end{tabular}}
		\end{table}
\subsection{Estimation of the parameters}
 Applying  the S$\&$P 500   set of samples and using the Bayesian MCMC method through Gibbs and griddy Gibbs sampling, we estimate the parameters    GARCH, STGARCH,   two-state MS-GARCH  and MS-STGARCH to compare their performance.  We also compare the proposed MS-STGARCH with EGARCH and GJR-GARCH. The prior density of transition probabilities  $\eta_{11}$ and $\eta_{22}$  are drawn from the beta distribution and  priors for the other parameters
 	  are assumed to be uniform over some finite intervals. We consider 10000 iterations  of Gibbs algorithm which half of them are burn-in-phase.   The posterior means and standard deviations for the parameters of the models corresponding to  S$\&$P 500 data  are reported in Table 6, which shows    that the  standard deviations are small enough in all cases, except for the slop parameter $\gamma_1$ which relates to the state with low volatility.

  The single-regime STGARCH has potential to react differently  to negative and positive shocks but
  does not consider shifting between
        different levels of volatilities.  The estimation results show that the level of volatility of the second regimes  in MS-STGARCH and MS-GARCH are higher. This is by the fact that the coefficient $a_{i2}$ for $j=0,1,2$ are respectively greater than $a_{i1}$.
                                 In Table 6,  we see that the estimated parameters of the MS-STGARCH satisfies  $a_{11}>a_{21}$ and $a_{12}>a_{22}$,  so the negative shocks have more affect on volatility than the positive shocks as $w_{j,t-1}$, for $j=1,2$ for large negative shocks approaches to zero

                                 Table 6 shows that   the posterior means of transition probabilities ($\eta_{11}$ and $\eta_{22}$) are close to one which indicate  less switch  between regimes. \textcolor{blue} {In Table 7 the estimated value of the parameters of EGARCH and GJR-GARCH are evaluated by applying the MLE method. }
                                    Estimated conditional transition  probabilities to the second state (high volatility regime), $\alpha_2^{(t)}$ computed by
                                    (3.6) plotted in Figure 4.

\label{fig3}
\input{epsf}
\epsfxsize=5.2in
\epsfysize=2in
 \begin{figure}
\centerline{\hspace{-.2in}\epsffile{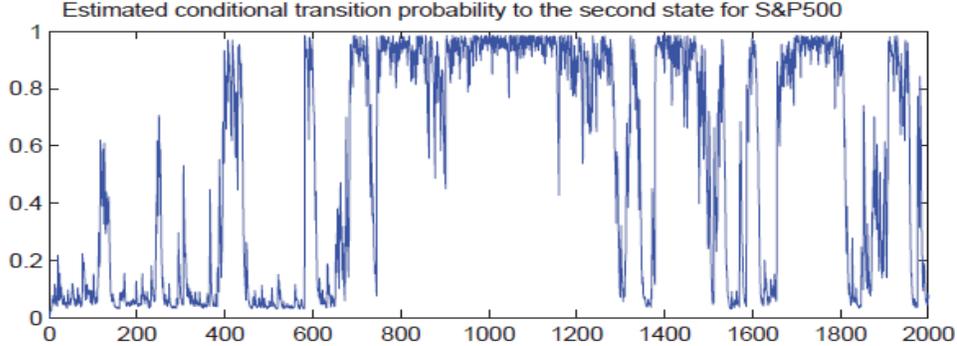}}
\vspace{0in}
\caption{\scriptsize Estimated conditional transition probabilities to the second states (the high volatility regime) of the fitted MS-STGARCH,  $\alpha_2^{(t)}$ calculated by (3.6), for the S$\&$P 500 log returns.}
\end{figure}
\begin{table}[!hbp]
	
	\caption{{\small Deviance information criterion (DIC)}}
	\hspace{3.5cm}\scriptsize
	\begin{tabular}{cc}
		\hline
		{\small Model} &{\scriptsize{ \small S$\&$P 500 returns }}\\
		\hline
		{\small GARCH} & 8464.8\\
		
		{\small ST-GARCH} & 7513.5\\
		{\small MS-GARCH} &7257.1\\
		{\small MS-STGARCH} &7147.8*\\
		\hline
	\end{tabular}
\end{table}

The MS-STGARCH has the potential to present better forecasting  when  different levels of volatilities are presented and  there is different effect for negative and positive shocks. The results of Table 8 demostrates that the MS-STGARCH has the best fitting to data. For appraising the performance of MS-STGARCH in forecasting, we survey the one-day-ahead value at risk (VaR) forecasts for the samples of S$\&$P 500.
    Based on the last 500 returns (of S$\&$P 500 ), the out
  of sample VaR \textcolor{blue}{forecasts} are calculated.
\label{fig4}
\input{epsf}
\epsfxsize=6.2in
\epsfysize=3.5in
 \begin{figure}
\vspace{0in}

\centerline{\hspace{0in}\epsffile{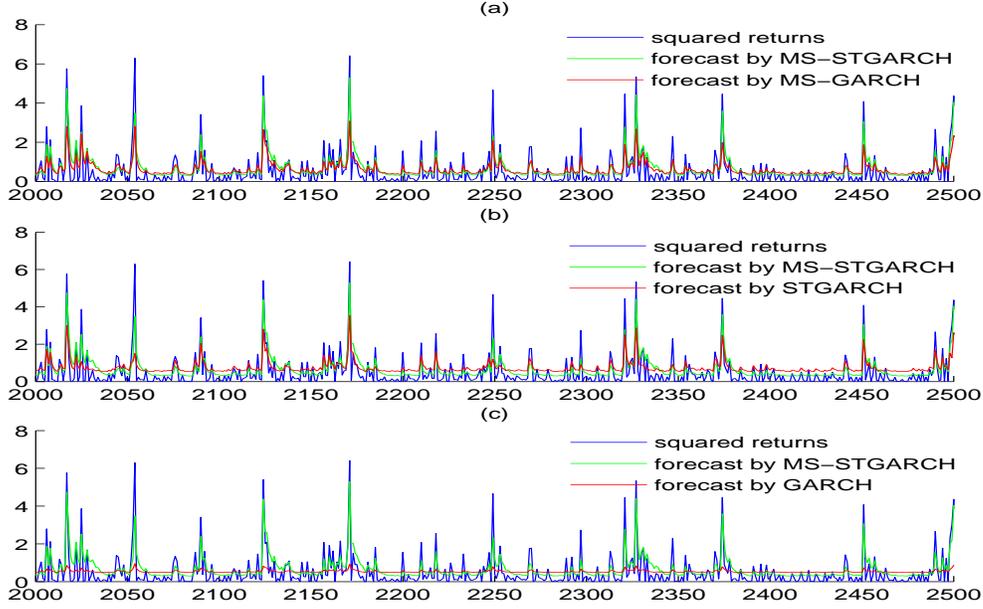}}
\vspace{0in}

\caption{\scriptsize Comparing  forecasts of MS-STGARCH  with forecasts of MS-GARCH, STGARCH and GARCH models  by applying squared returns  of  S$\&$P 500 (blue)  (a):  forecast by MS-STGARCH (green) and forecast by MS-GARCH (red), (b) forecast by MS-STGARCH (green) and forecast by STGARCH (red), (c)  forecast by MS-STGARCH (green) and forecast by GARCH (red).}
\end{figure}


\begin{table}
	\vspace{-.05in}
	\caption{{\scriptsize VaR results  of  S$\&$P 500  daily log returns.}}
	\hspace{-.5cm}{\hspace{.8in}\scriptsize
		\begin{tabular}{c c c c c c c}
			\hline
			Model   & $\alpha$ & $E(V_t(\alpha))$ & N & UC & IND &  CC  \\
\hline
                    &0.99 &5 &9&2.596  & 0.330   & 2.926 \\
			&0.95 &25 &26&0.038 &2.740 &2.778\\
			&0.9 &50 &42 &1.531 &0.070 &1.601\\
			EGARCH   &0.1&50&32&8.227 &0.001&8.228\\
			&0.05&25&12&8.790&1.155&9.946\\
			&0.01&5&1&4.829& 0.004 & 4.833 \\

                    \hline
                   &0.99 &5 &10&3.891  & 0.408   & 4.299 \\
			&0.95 &25 &26&0.037 &2.740 &2.778\\
			&0.9 &50 &42 &1.531 &0.070 &1.601\\
			GJR-GARCH   &0.1&50&31&9.236 &0.036&9.239\\
			&0.05&25&16&3.925&2.775&6.700\\
			&0.01&5&3&0.950& 0.036 & 0.987 \\

			\hline
			&0.99 &5 &9&2.596  & NA   & NA \\
			&0.95 &25 &18 &2.276 &2.024 &4.3\\
			&0.9 &50 &44 &0.830 &0.252 &1.082\\
			GARCH   &0.1&50&42&1.395 &0.071&1.466\\
			&0.05&25&27&2.276&0.142&2.418\\
			&0.01&5&4&0.229& NA & NA \\
			\hline
			&0.99&5&5&0& NA & NA \\
			&0.95&25&20&1.147& NA & NA \\
			&0.9&50&31&9.235&0.5317&9.767\\
			STGARCH&0.1&50&28&12.684&1.191&13.875\\
			&0.05&25&8&16.441& NA & NA \\
			&0.01&5&2&2.365& NA & NA \\
			\hline
			&0.99&5&8&1.526& NA & NA \\
			&0.95&25&27&0.156&0.1787&0.335\\
			&0.9&50&42&1.531&0.827&2.358\\
			MS-GARCH&0.1&50&37&4.112&0.252&4.364\\
			&0.05&25&15&4.926&0.539&5.465\\
			&0.01&5&2&2.365& NA & NA \\
			\hline
			&0.99 &5 &9&2.596  & NA  & NA \\
			
			&0.95 &25 &27 &0.156&0.1787&0.335\\
			&0.9 &50 &45 &0.595 &0.301 &0.897\\
			MS-STGARCH   &0.1&50&44&0.9&0.252&1.109\\
			&0.05&25&18&2.3&0.275&2.575\\
			&0.01&5&2&2.38& NA & NA \\
			\hline

		\end{tabular}}
	\end{table}

\par   According to the results of Table 9, at the $ 5\%$ and $10\%$ significance levels, the  $LR_{uc}$ test is
	rejected three times  for EGARCH, GJR-GARCH and  STGARCH, two times for MS-GARCH   and is accepted  at all risk levels  $\alpha$  for the GARCH and MS-STGARCH models.
  The $LR_{ind}$ statistic at $ 5\%$ significance level is smaller than critical  value 
$\chi^2_{0.95}$ with one degree of freedom  for all cases that test is applied. Also excluding the cases of risk level 0.95 for the  EGARCH and GJR-GARCH, the IND test is accepted at  $10\%$. At the $ 5\%$ significance level, the conditional coverage (CC) test is higher than  critical value $\chi^2_{0.95}$ with two degrees of freedom two times for the EGARCH, GJR-GARCH and STGARCH while  at the $ 10\%$ significance level this test is rejected three times for the  EGARCH, two times for the GJR-GARCH and STGARCH and one time for the MS-GARCH.

\par 
In order to appraise the ability of competing models to forecast volatility, we apply the Diebold Mariano test. Testing for equal forecast accuracy is an approach to evaluate the predictive capability of competitor models. Diebold and Mariano (1995) proposed a unified method for testing the null hypothesis of no difference in the forecasting accuracy of two competing models \cite{xekalaki}.  Harvey, Leybourne and Newbold (1997) suggested a modified of Diebold Mariano (DM) test for small sample.  The DM test or its HLN variant are applied widely in empirical forecasting research, see \cite{deibold}, \cite{harvey} and \cite{curto b}. Consider two forecast sequences as
	\begin{equation*}
	\{\hat{y}_{it}: t=1,\cdots,T\},  i=1,2;
	\end{equation*}
		and define
		\begin{equation*}
		e_{it}=\hat{y}_{it}-y_t
		\end{equation*}
	that $\{y_t,t=1,\cdots,T\}$ are actual values. Let $g(e_{it})=e^2_{it}$ and 
	\begin{equation*}
	d_t=g(e_{1t})-g(e_{2t});
	\end{equation*}
we would like to test the null hypothesis:
\begin{equation*}
H_0: E(d_t)=0,  \forall t
\end{equation*}			
versus the alternative hypothesis
\begin{equation*}
H_1: E(d_t)<0,  
\end{equation*},
under covariance stationarity of the process $\{d_t:1,\cdots,T\}$, the Diebold-Mariano (DM) statistic for testing the null hypothesis is given by:
\begin{equation*}
\frac{\bar{d}}{\sqrt{\hat{Var}((\bar{d})}},
\end{equation*}
and is approximately normally distributed for large samples. 
For evaluating the performance of our model in one-step ahead conditional variance forecast, we compute the DM statistic for pairwise comparison of MS-STGARCH model  with GARCH, STGARCH, MS-GARCH,GJR and EGARCH models. According to the test  results demonstrated  in table 10, the null hypothesis  is rejected at the $5 \%$ significance level for all cases as all the statistics are less than $Z_{0.05}$. So our presented model has an improvement in the forecasting performance.
\begin{table}
	\caption{{\scriptsize DM test results}}
	\hspace{1.5cm}{\hspace{.8in}\scriptsize
		\begin{tabular}{cc }
			\hline
		Comparison of MS-STGARCH with  &  Statistic value \\
					\hline
                                EGARCH & -3.89\\
                               GJR-GARCH & -3.73\\
					GARCH & -3.84\\
					STGARCH & -4.08\\
					MS-GARCH &  -2.08\\
					\hline
						\end{tabular}}
	\end{table}}
Also for specifying  the out of sample forecast performance of  the MS-STGARCH toward the competing models, We  compare the forecasting volatility  $E(Y^2_t|\mathcal{F}_{t-1})$, or conditional variance, of GARCH, STGARCh and MS-GARCH   with the squared returns. In Figure 5, the squared returns of $S\&P500$ and the   forecasting values of  competing models are plotted. According to this figure the differences of forecast  and real values (errors) in the MS-STGARCH always has been much less than other compared models.  The results of  Table 11 show that the least values of the MSE and MAE are related to   the MS-STGARCH model that  reveals   the best  forecast compared with the  other reviewed  models in this paper.

\begin{table}[!hbp]

\caption{{\small Measures of performance forecasting}}
\hspace{1.5cm}\scriptsize
\begin{tabular}{ccc}
  \hline
 {\small Model} &{\scriptsize{ \small Mean square error (MSE) }}& Mean absolute error (MAE)\\
  \hline

{\small EGARCH} &0.676  &0.528\\
  {\small GJR-GARCH} &0.665  &0.517\\
{\small GARCH} &0.707  & 0.524 \\
   {\small ST-GARCH} &0.384 &0.498 \\
 {\small MS-GARCH} &0.302 & 0.395\\
 {\small MS-STGARCH} & 0.227*  & 0.369*\\
  \hline
\end{tabular}
\end{table}
\par 







\section{Conclusion}
\textcolor{blue}{ Applying Markov switch structure cause to have a better fitting  while the existence of different levels of volatilities are evident.  
 Also the asymmetry  effects of negative and positive shocks in many case are trivial and transition between this effects happens in some smooth ways and not sudden. So in many cases the  use of MS-STGARCH has advantages to the other methods as GARCH, MS-GARCH, ST-GARCH, EGARCH and GJR-GARCH as we find this performance for SP500 indicies where studied in the paper.
In such cases  a much better fit to the data can be provided by the presented model which leads to  the forecasts 
with much smaller  error.}
The MS-STGARCH extends the MS-GARCH model by considering convex combination of time dependent logistic weights between the effect of negative and positive shocks  in each regime. It  also extends  the STGARCH model by considering  transition between different levels of volatilities.   We show that the existence of  a simple  condition  suffices for the existence of an 
asymptotic upper bound for the second moment of returns which causes the stability of the model.

 By fitting  the GARCH, STGARCH, MS-GARCH and MS-STGARCH models to the   S$\&$P 500 log returns 
 we find that our model has the best  DIC, see Table 8, and provides the best forecast volatilities, see Figure 5.  Also  in performing  Diebold Mariano test, our model has  the best performance in compare to 
 the 
EGARCH and GJR-GARCH. 
The forecasted one-day ahead Value at Risk (VaR) of our model has better performance to the EGARCH,GJR-GARCH, STGARCH and MS-GARCH.

\par  Further researches could be oriented to investigate the  existence of the  third and the fourth moments of the process and derive the necessary and sufficient conditions for  stationarity and  ergodicity  of the process.  For the sake of simplicity it was assumed that the process conditional mean is zero, this assumption could be relaxed by refining the structure of model to allow ARMA structure for conditional mean. 
Since Financial time series data are typically observed to have heavy tails \cite{liu}, it might be interesting to replace the Gaussian distribution with Student's t  or
stable Paretian distributions to investigate the ability for  modeling heavy  tailed  property  of financial time series such as \cite{curto a} work.
\newpage
\begin{center}
{\bf Appendix A}
\end{center}
Proof of Remark 3.1.\\
\\
As the hidden variables $\{Z_{t}\}_{t\geq 1}$ have Markov structure in MS-STARCH model, so
\begin{align}
\alpha_j^{(t)}= &p(Z_t=j|\mathcal{I}_{t-1})=\sum_{m=1}^{K}{P(Z_t=j,Z_{t-1}=m|\mathcal{I}_{t-1})}
\nonumber\\
&=\sum_{m=1}^{K}{p(Z_t=j|Z_{t-1}=m,\mathcal{I}_{t-1})p(Z_{t-1}=m|\mathcal{I}_{t-1})}
\nonumber\\
&=\sum_{m=1}^{K}{p(Z_t=j|Z_{t-1}=m)p(Z_{t-1}=m|\mathcal{I}_{t-1})}
\nonumber\\
&=\frac{\sum_{m=1}^{K}{f(\mathcal{I}_{t-1},Z_{t-1}=m)p_{m,j}}}{\sum_{m=1}^{K}{f(\mathcal{I}_{t-1},Z_{t-1}=m)}}
\nonumber\\
&=\frac{\sum_{m=1}^{K}{f(y_{t-1}|Z_{t-1}=m,\mathcal{I}_{t-2})}\alpha_m^{(t-1)}p_{m,j}}{\sum_{m=1}^{K}{f(y_{t-1}|Z_{t-1}=m,\mathcal{I}_{t-2})\alpha_m^{(t-1)}}}.
\end{align}
\begin{center}
{\bf Appendix B}
\end{center}
Proof of Theorem 3.1.\\
Let $E_{t}(.)$  denotes the expectation with respect to  the information up to time t. Thus the second moment  of the model can be calculated as, \textcolor{blue}{see \cite{abramson}}:
$$E(y^2_t)=E(H_{Z_t,t})=E_{Z_t}[E_{t-1}(H_{Z_t,t}|z_t)]\hspace{6cm}$$

\begin{equation}\label{10}
=\sum_{z_t=1}^{K}{\pi_{z_t}E_{t-1}(H_{Z_t,t}|z_t)}.\hspace{2.5cm}
\end{equation}

Also
 let  $E(.|z_t)$ and $p(.|z_t)$ denote $E(.|Z_t=z_t)$ and $P(.|Z_t=z_t)$, respectively, where $z_t$ is the realization of the state at time t.
 \textcolor{blue}{Applying to the method of Medeiros      \cite{medeiros}}, we find an upper bound of $E_{t-1}(H_{m,t}|z_t)$, for  $ m=1,2,\cdots ,K$ by the following
\begin{align}\label{11}
E_{t-1}(H_{m,t}|z_t)=
&E_{t-1}(a_{0m}+a_{1m}y^2_{t-1}(1-w_{m,t-1})
+a_{2m}y^2_{t-1}w_{m,t-1}+b_mH_{m,t-1}|z_t)
\nonumber\\
= &\underbrace{a_{0m}}_{I}+\underbrace{a_{1m}E_{t-1}[y^2_{t-1}|z_t]}_{II}+\underbrace{(a_{2m}-a_{1m})E_{t-1}[y^2_{t-1}w_{m,t-1}|z_t]}_{III}
\nonumber \\
&\qquad\qquad\qquad \ \ \ \ \ \ \ \ \qquad\qquad+\underbrace{b_mE_{t-1}[H_{m,t-1}|z_t]}_{IV}.
\end{align}
The term (II) in (\ref{11}) can be interpreted  as follows:

 $$E_{t-1}[y^2_{t-1}|z_t]=\sum_{z_{t-1}=1}^{K}
 {\int_{S_{\mathcal{I}_{t-1}}}{y^2_{t-1}p(\mathcal{I}_{t-1}|z_t,z_{t-1})p(z_{t-1}|z_t)d\mathcal{I}_{t-1}}}$$
 \begin{equation}\label{12}
 =\sum_{z_{t-1}=1}^{K}{p(z_{t-1}|z_t)E_{t-2}[H_{Z_{t-1},t-1}|z_{t-1}]},
 \end{equation}
 where $S_{\mathcal{I}_{t-1}}$ is the support of $\mathcal{I}_{t-1}=(y_1,\cdots,y_{t-1})$.
 \\
{\bf Upper bound for III in  (\ref{11}):}
Let  $0<M<\infty$ be a constant, so
\begin{align}
E_{t-1}[y^2_{t-1}w_{m,t-1}|z_t]= &E_{t-1}[y^2_{t-1}w_{m,t-1}I_{|y_{t-1}|<M}|z_t]
\nonumber\\
&+E_{t-1}[y^2_{t-1}w_{m,t-1}I_{|y_{t-1}|\geq M}|z_t]
\nonumber
\end{align}
in which
\begin{displaymath}
I_{x<a}=\left\{ \begin{array}{ll}
1 & \textrm{if  $x<a$}\\
0 & \textrm{otherwise.}
\end{array} \right.
\end{displaymath}
As by (\ref{4}),  $0<w_{m,t-1}<1$ and so
\begin{equation*}
E_{t-1}[y^2_{t-1}w_{m,t-1}|z_t]\leq M^2+E_{t-1}[y^2_{t-1}w_{m,t-1}I_{|y_{t-1}|\geq M}|z_t],
\end{equation*}
also
\begin{align}
E_{t-1}[y^2_{t-1}w_{m,t-1}I_{|y_{t-1}|\geq M}|z_t]=&\int_{S_{\mathcal{I}_{t-2}},y_{t-1}\leq -M}{y^2_{t-1}[w_{m,t-1}]p(\mathcal{I}_{t-1}|z_t)d\mathcal{I}_{t-1}}
\nonumber\\
&+\int_{S_{\mathcal{I}_{t-2}},y_{t-1}\geq M}{y^2_{t-1}[w_{m,t-1}]p(\mathcal{I}_{t-1}|z_t)d\mathcal{I}_{t-1}},
\nonumber
\end{align}
by (\ref{4}),
\begin{equation}
\lim_{y_{t-1}\rightarrow +\infty}w_{m,t-1}=1 \ \ \ \ \and\ \ \ \ \lim_{y_{t-1}\rightarrow -\infty}w_{m,t-1}=0.
\end{equation}
So for any fixed positive  small number $\delta>0$, we can consider  $M>0$ so large that for  $y_{t-1}\geq M$, $|w_{m,t-1}-1|\leq \delta$ and for $y_{t-1}\leq -M$, $|w_{m,t-1}|\leq \delta$. Hence
\begin{align}
E_{t-1}[y^2_{t-1}w_{m,t-1}I_{|y_{t-1}|\geq M}|z_t] \leq \ & \delta \int_{S_{\mathcal{I}_{t-2}},y_{t-1}\leq -M}{y^2_{t-1}p(\mathcal{I}_{t-1}|z_t)d\mathcal{I}_{t-1}}
\nonumber\\
&+ (\delta+1) \int_{S_{\mathcal{I}_{t-2}},y_{t-1}\geq M}{y^2_{t-1}p(\mathcal{I}_{t-1}|z_t)d\mathcal{I}_{t-1}}.
\nonumber
\end{align}
Since the distribution of the $\{\varepsilon_{t}\}$ is symmetric, then
\begin{align}
\delta\int_{S_{\mathcal{I}_{t-2}},y_{t-1}\leq -M}{y^2_{t-1} p(\mathcal{I}_{t-1}|z_t)d\mathcal{I}_{t-1}}\leq & \delta \int_{S_{\mathcal{I}_{t-2}},-\infty<y_{t-1}<0}{y^2_{t-1}p(\mathcal{I}_{t-1}|z_t)d\mathcal{I}_{t-1}}
\nonumber\\
&=\delta\frac{ E_{t-1}[y^2_{t-1}|z_t]}{2}
\nonumber
\end{align}
\\
and
\begin{align}
(\delta+1)\int_{S_{\mathcal{I}_{t-2}},y_{t-1}\geq M}{y^2_{t-1}p(\mathcal{I}_{t-1}|z_t)d\mathcal{I}_{t-1}}\leq & (\delta+1)\int_{S_{\mathcal{I}_{t-2}},0<y_{t-1}<\infty}{y^2_{t-1} p(\mathcal{I}_{t-1}|z_t)d\mathcal{I}_{t-1}}
 \nonumber\\
 &=(\delta+1)\frac{ E_{t-1}[y^2_{t-1}|z_t]}{2}.
 \nonumber
\end{align}
\\
Therefor
\begin{equation*}
E_{t-1}[y^2_{t-1}w_{m,t-1}|z_t]\leq  M^2+(\delta+\frac{1}{2})E_{t-1}[y^2_{t-1}|z_t].
\end{equation*}
{\bf Upper bound for IV in  (\ref{11}):}
\begin{equation*}
  b_{m}E_{t-1}(H_{m,t-1}|z_t)= b_{m}{\int_{S_{\mathcal{I}_{t-1}}}{H_{m,t-1}p(\mathcal{I}_{t-1}|z_t)d\mathcal{I}_{t-1}}}
    \end{equation*}
    \begin{equation}
    = b_{m}\sum_{z_{t-1}=1}^{K}
 {p(z_{t-1}|z_t)E_{t-2}(H_{m,t-1}|z_{t-1})}.
 \end{equation}
  By replacing the obtained  upper bounds and relations (\ref{12}) in (\ref{11}),  the upper bound for $E_{t-1}(H_{m,t}|z_t)$ is obtained  by:
\begin{align}
E_{t-1}(H_{m,t}|z_t)&\leq a_{0m}+|a_{2m}-a_{1m}|M^2
\nonumber\\
&+\sum_{z_{t-1}=1}^{K}{{[a_{1m}+|a_{2m}-a_{1m}|(\delta+\frac{1}{2})]
p(z_{t-1}|z_t)}E_{t-2}[H_{z_{t-1},t-1}|z_{t-1}]}\nonumber
\\&+\sum_{z_{t-1}=1}^{K}{b_mp(z_{t-1}|z_t)E_{t-2}(H_{m,t-1}|z_{t-1})},
\end{align}
in which by Bayes' rule
$$p(z_{t-i}|z_t)=\frac{\pi_{z_{t-i}}}{\pi_{z_t}}\{P_{z_{t-i}z_t }\},$$
where $P$ is the transition probability matrix.

Let $A_t(j,k)=E_{t-1}[H_{j,t}|Z_{t}=k]$, $ {A}_t=[A_t(1,1),A_t(2,1),\cdots,A_t(K,1),A_t(1,2),\\\cdots,A_t(K,K)]$ be a $K^2$-by-1  vector and consider $\dot{ { {\Omega}}}=( { {\Omega}}^{\prime},\cdots, { {\Omega}}^{\prime})^{\prime}$ be a  vector that is made of K vector  ${ {\Omega}}$.
By  (7.23)-(7.26), the following recursive inequality is attained,
\begin{equation}\label{19}
\textbf{A}_t\leq \dot{ { {\Omega}}}+ {\bf{C}}\textbf{A}_{t-1},\ \ \ t\geq 0.
\end{equation}
with some initial conditions $\textbf{A}_{-1}.$
{
{  The relation (\ref{19}) implies that
\begin{equation}\label{20}
 {A}_t\leq \dot{{ {\Omega}}}\sum_{i=0}^{t-1}{{ {C}}^i}+{ {C}}^{t} {A}_{0}:= {B}_{t}.\hspace{8.5cm}
\end{equation}
Following the matrix convergence theorem \cite{lancaster}, the necessary  condition for the convergence of $ {B}_{t}$ when $t\rightarrow\infty$
is that $\rho({ {C}})
<1$. Under this condition, ${ {C}}^t$ converges to zero as t goes to infinity and $\sum_{i=0}^{t-1}{{ {C}}^i}$ converges to $(I-{ {C}})^{-1}$
provided that matrix $(I-{ {C}})$ is invertible.   So if $\rho({ {C}})
<1$,
$$\lim_{t\rightarrow\infty} {A}_t\leq(I-{ {C}})^{-1}\dot{{ {\Omega}}}.\hspace{8.5cm}$$
  By (\ref{10})
the upper bound for the asymptotic behavior of  unconditional variance is given by
$$lim_{t\rightarrow\infty}E(y^2_t)\leq{ {\Pi^\prime(I-C)^{-1}}}\dot{{ {\Omega}}}.\hspace{8cm}$$}}\\\

\noindent
{{\bf Acknowledgement:} This paper initiated  during Professor Rezakhah sabbatical at Institute of Mathematics at the EPFL where benefitted from the discussion of the paper with Professor Stephan Morgenthaler and
 careful written comments and suggestions of  Professor  Anthony Davison   which  caused to improve the quality of this paper. }

\bibliographystyle{plain}

\end{document}